\input amstex
\input xy
\xyoption{all}
\documentstyle{amsppt}
\document
\magnification=1200
\NoBlackBoxes
\nologo
\hoffset1.5cm
\voffset2cm
\vsize15.5cm
\def\F{\bold{F}}

\def\C{\bold{C}}
\def\P{\bold{P}}
\def\Q{\bold{Q}}

\def\Z{\bold{Z}}
\def\R{\bold{R}}
\def\F{\bold{F}}
\def\N{\bold{N}}

\def\cM{\Cal{M}}

\def\cO{\Cal{O}}

\bigskip

%\hfill{\it draft 15.09.2020}
%\hfill{\it draft 26.12.2020}

%\hfill{\it BirationalAndNori.tex}
\bigskip

\centerline{\bf BIRATIONAL MAPS AND NORI MOTIVES}

\medskip

\centerline{\bf No\'emie~C.~Combe,\quad  Yuri~I.~Manin,\quad  Matilde~Marcolli}

\bigskip

 ABSTRACT. The monograph [HuM-St17]  contains a systematical exposition of  Nori motives 
 that were developed and studied as the ``universal (co)homology
 theory'' of algebraic varieties (or schemes), according to the prophetic vision
 of A. Grothendieck. Since then, some research
 was dedicated to application of Nori motives in various domains
 of algebraic geometry: geometries in characteristic 1 ([LieMaMar19], [MaMar18]),
 absolute Galois group ([MaMar19-2]), persistence formalism ([MaMar19-1]).
 
 In this note, we sketch an approach to the problems
of equivariant  birational geometry developed by M. Kontsevich and Yu. Tschinkel  in [KTsch20],
where Burnside invariants were introduced. We are
making explicit the role of Nori constructions in this environment.


\bigskip

\centerline{\bf  0. INTRODUCTION AND SUMMARY}

\medskip

{\bf 0.1. Birational maps and their symmetries.} Our main objects of study here are stable {\it birational maps}, mostly between
algebraic varieties defined over a subfield of $\C$.
The general restriction of {\it stability} is discussed in [AVi02], and main results
about natural categories/``towers''  of birational maps we use here
are given in [AT16] and [BeRy19].
\smallskip

{\it Symmetries} of stable maps and their moduli spaces appear in various contexts.
The celebrated Grothendieck's approach to the study of {\it absolute Galois group} $G_{\Q}$
of the field of all algebraic numbers $\overline{\Q}$ bridged geometry and arithmetic via the
tower of stable \'etale maps to $\P^1\setminus \{0,1, \infty \}$. More generally,
for any integral scheme $X$, the exact sequence
$$
1 \to \pi_1(X\otimes_{\Q}\overline{\Q}) \to  \pi_1(X) \to G_{\Q} \to 1
\eqno(0.1)
$$
merges actions of what we call  {\it arithmetic} and {\it geometric} symmetries.    
\smallskip

In [CMMar20], in place of \'etale maps $\P^1\setminus \{0,1, \infty \}$ we considered
moduli spaces of stable curves of genus zero with marked points, and later demonstrated that
their geometric symmetries can be approached through the {\it Nori motivic structures}
of these moduli spaces. In this context, the exact sequence of fundamental groups above
is replaced by a subextension of the ``motivic fundamental groups''
$$
1\to G_{mot}(\overline{\Q} ,\Q) \to G_{mot}(\Q ,\Q) \to G_{\Q} \to 1
\eqno(0.2)
$$
Here we focus on similar constructions, but starting with towers of birational maps
replacing towers of stable moduli spaces above.

\bigskip

{\bf 0.2. Burnside groups and Nori diagrams.} 
In [AT16], [BeRy19], and other papers it was shown that the problem 
of classification of birational maps can use information encoded
in the tower of natural maps structurally similar to the tower 
of stable genus zero modular spaces $\overline{M}_{0,n}$. Using these natural maps,
in the recent article [KrTsch20], A. Kresch and Yu. Tschinkel are 
 imposing an additional geometric symmetry group from the start and then showing that
this geometric symmetry can be encoded by certain analogues of ``modular symbols''.

\smallskip

In this note we demonstrate that the Kresch--Tschinkel
modular symbols defined through {\it Burnside  groups},
also have a natural description in terms of Nori motivic
structures. Briefly, we prove the following theorem:

\smallskip

{\bf 0.2.1. Theorem.} {\it Nori stratifications of appropriate towers of birational maps
can be enriched to the homological Nori geometric diagram by
Burnside groups of strata.}

\smallskip

For a more precise statement, see Theorem 3.3.

\bigskip

{\bf 0.3. Quantum statistical mechanics via Burnside groups.} Finally,
in the last section we consider the Bost--Connes formalism connecting
arithmetic zeta--functions of fields of algebraic numbers
with physicists' studies of classical and quantum behaviour
of physical systems with infinite number of degrees of freedom:
see [BoCon95];  Sec. 5 of [MaMar18], and references therein.

\smallskip

It was In this Sec. 5, [MaMar18], that  we have sketched the enrichment of
Bost--Connes via Burnside groups which  in more detail will be presented here.

\bigskip

{\bf 0.4. Summary.} Many preliminaries and necessary definitions in this paper are distributed
as follows. Sec. 1 contains a short dictionary of graphic presentations of categories and functors.
 Nori diagrams and representations are discussed in Sec. 2;
Burnside groups in Sec. 3, containing also the last steps of the proof of Theorem 0.2.1.
Section 4 recalls a construction of [LieMaMar19] of a lifting of the
Bost--Connes algebra to the relative $\hat\Z$-equivariant Kontsevich-Tschinkel Burnside
group and to its underlying assembler category and spectrum, and to equivariant
Nori motives. Motivated by this construction, we then show that a similar
Bost-Connes-type structure is also present in the Kontsevich-Pestun-Tschinkel
modular symbols. 

\bigskip
\centerline{\bf 1. BACKGROUND}

\bigskip

{\bf 1.1. Categories and their diagrams.}  A diagram $D$ is family $(V(D),E(D), \partial )$
where $\partial $ (boundary map) is an embedding of $E(D)$ (edges) into
$V(D)\times V(D)$ (ordered pairs of vertices). Each category defines its diagram, whose 
edges are its morphisms, and vertices are its objects.

\smallskip

Conversely, diagrams themselves form objects of category, whose morphisms
imitate functors.

\smallskip

For a more detailed discussion of combinatorics of categories
based upon  diagrams, we refer the reader to 
Sec. 0.2 of [MaMar19-1] and Sec. 5 of [CM19].  Here we remind only the definition
of {\it posets in groupoids} ([CM19], Def. 5.2).
\medskip

{\bf 1.2. Definition.} {\it A  category $\Cal{P}\Cal{G}$ is called
a poset in groupoids, if 

\smallskip

(a) For any object $X$ of  $\Cal{P}\Cal{G}$ , the full subcategory consisting
of all objects isomorphic $X$, is a groupoid, that is, all morphisms in it
are isomorphisms.

\smallskip

(b) Whenever $X_1$ and $X_2$ are  not isomorphic and $\roman{Hom} (X_1, X_2) \ne \emptyset$,
then $\roman{Hom} (X_1, X_2)$ has a single orbit with respect to the left action of the
group  
$$
\roman{Hom} (X_1, X_1)\times \roman{Hom} (X_2, X_2)^{op}
$$
combining precomposition and postcomposition.}

\smallskip

A part of posets in groupoids consists of {\it thin categories} $\Cal{C}$: such that 
if any set $\roman{Hom}_{\Cal{C}} (X_1,X_2)$ has cardinality $\le 1$, and if both
$\roman{Hom}_{\Cal{C}} (X_1,X_2)$ and $\roman{Hom}_{\Cal{C}} (X_2,X_1)$
are non--empty, then $X_1 =X_2$.

\smallskip

\medskip

{\bf 1.3.  Diagrams of effective pairs.} ([HuM-St17], Ch. 9, Def. 9.1.1.) Fix a subfield $k$ of $\bold{C}$
and define the diagram $\roman{Pairs^{eff}}$ of effective pairs over $k$ in the following way.
\smallskip

(a) Vertices of  $\roman{Pairs^{eff}}$ are triples $(X,Y,i)$ where $X$ is a variety over $k$,
$Y\subset X$ is a closed subvariety, and $i\in \Z$.

\smallskip

(b) There are two types of edges of  $\roman{Pairs^{eff}}$: {\it functoriality edges} and
{\it coboundary edges}:

\smallskip

(b1) Each morphism $f: X\to X^{\prime}$ with $f(Y)\subset Y^{\prime}$ determines  edges
denoted $(f^*,i)$ starting at $(X,Y,i)$ and landing at $(X^{\prime}, Y^{\prime}, i+1)$ for
every $i\in \Z$.

\smallskip

(b2) Each ladder $Z\subset Y\subset X$ of closed subvarieties determines edges $(\partial ,i)$
starting at $(Y,Z,i)$ and ending at $(X,Y, i+1)$ for every $i \in \Z$.

\smallskip

Later, whenever we will have to consider tensor structures on the categories of diagrams,
we will have to consider (super)gradings of such diagrams which are discussed in [HuM-St17],
Ch. 8.

\medskip

{\bf 1.4.  Categories of blowings up.} We will now describe some categories of good blowings up,
following Sec. 1.3 of [AT16] (with somewhat changed terminology and notation).

\smallskip

 We will call a ``good'' scheme
what in the Introduction to  [AT16] is called a  ``noetherian quasi excellent (qe) regular scheme''.

Consider a morphism of good
schemes $\varphi : X_1 \to X_2$, which is the blowing up of a coherent sheaf of ideals $I\subset \cO_{X_2}$
(for the relevant definitions in this context, see [AT16], Sec. 2.1.8). Alternatively,
we will call such $\varphi$ the blowing up of the closed subscheme defined by equations the $f=0$
for all $f\in I$.

\smallskip

Assume also given normal crossings divisors $D_i\subset X_i$ such that $D_1 = \varphi^{-1} (D_2)$.

\smallskip

Let $U$ be the maximal open subscheme of $X_2$ upon which the restriction
of $I$ is its structure sheaf.
%and whose boundary is disjoint from $U$. 
It follows
that $\varphi$ induces an isomorphism $\varphi^{-1} (U) \to U$.

\smallskip

We will call {\it a good morphism} the structure represented by a set of data $(X_i,D_i,I, \varphi )$
as above. In particular, identical morphisms, and generally, automorphisms, are good:
for them $D_1=D_2=\emptyset$. 

\medskip

Consider a finite connected poset in groupoids $\Cal{M}$ whose objects are
data $(X,D,I)$ and
morphisms are good morphisms $\varphi$  in the sense of Sec. 0.1 above.

%\smallskip
%Assume also for simplicity that all $X$ in the data are smooth projective manifolds.

\smallskip

Consider the diagram (with identities) $D(\Cal{M})$ whose vertices are objects
of $\Cal{M}$, oriented edges are morphisms of $\Cal{M}$, and orientation of
$X\to Y$ is from $X$ to $Y$.
Call edges corresponding to isomorphisms (in particular, identities) {\it horizontal} ones, and 
other edges {\it vertical } ones. 

\smallskip

Form also the following quotient of $D(\Cal{M})$ which we denote $T(\Cal{M})$:
vertices of  $T(\Cal{M})$ are isomorphisms classes of $\Cal{M}$, and 
oriented edges of $T(\Cal{M})$ are orbits of non--empty sets $\roman{Hom}_{\Cal{M}}(X_1,X_2)$
with respect to the left action of the
group  
$$
\roman{Hom} (X_1, X_1)\times \roman{Hom} (X_2, X_2)^{op}
$$
combining precomposition and postcomposition,
as in Def. 1.2 above. We omit edges corresponding to identities.
\smallskip

Call an edge $X_1\to X_2$ in $T(\Cal{M})$ {\it indecomposable},
if the respective morphism cannot be expressed as composition of
other morphisms.

\smallskip

In an important particular case, the diagram $T(\Cal{M})$ is in fact a tree  oriented downwards. More precisely,
starting with any its vertex $X_0$, we may consider the longest
sequence of vertices (``a path down'')
$$
X_0 \to X_1\to \dots \to  X_h.
$$
Assume that there is only one vertex from which any longest path
down can start.

\bigskip

{\bf 1.4.1. Definition.} One object of the category $\roman{Bl}_{\roman{rs}}$  ({\it the regular surjective category of
blowings up}, cf. Sec. 1.3 of  [TM16]) is a triple $(X_2, I, D_2)$ constituting
a part of good morphism $\varphi$ as above.
\smallskip
One morphism  between such objects of $\roman{Bl}_{\roman{rs}}$
$$
(X_i^{\prime},D_i^{\prime} ,I^{\prime}, \varphi^{\prime}) \to  (X_i,D_i,I, \varphi )
$$
is represented by a regular and surjective morphism
$$
g: X_2^{\prime}\to X_2
$$
satisfying the following conditions:
$$
g^{-1}(D_2)= D_2^{\prime}, \quad\quad   g^*(I) = I^{\prime}. 
$$

\medskip

{\it Remark.} From this definition one can deduce, that $g$ induces a canonical isomorphism
$$
X_1^{\prime} \to X_1 \times_{X_2} X_2^{\prime},
$$
and moreover, $D_1^{\prime}$ is the inverse image of $D_1$ with respect to the
composition of this isomorphism and projection $X_1 \times_{X_2} X_2^{\prime}\to X_1$
(cf. Definition 1.3.1 of [AT16]). This presentation might me helpful for defining
and studying compositions of morphisms between objects of $\roman{Bl}_{\roman{rs}}$.

\smallskip

In the next Sections, we will  be studying posets of groupoids as above
from the viewpoint of Nori theory, as presented in Sec.1 of [MaMar19-1].
But before starting it, we need one more definition.

\medskip

{\bf 1.5. Simple normal crossings divisors.} Let $S$ be a finite set. Consider a family
 $\{D_s \subset X\,|\,s\in S\}$ of closed immersions. Following Definition (3.1) of [BeRy19],
 we will call it an {\it $S$--labelled simple normal crossings divisor on $X$}, if
 for any finite subset $S^{\prime}\subseteq S$ the intersection $\bigcap_{s\in S^{\prime}}D_s$
is smooth of codimension $\roman{card}\, S^{\prime}$.
\bigskip 

\centerline{\bf 2. NORI GEOMETRIC DIAGRAMS OF BLOWINGS UP}

\bigskip

{\bf 2.1. Categories $\Cal{M}$.} We will now introduce a class of geometric categories
that will be the starting point for our enrichment  of birational maps by Nori motives.
Notation $\Cal{M}$ for a generic member of this class should remind the reader
that we generalise here the basic example of stable modular spaces of genus zero
and their canonical stratifications studied in [CMMar20]. Here are the basic restrictions
imposed upon $\Cal{M}$.

\smallskip

(a) Objects of $\Cal{M}$ are some objects of  $\roman{Bl}_{\roman{rs}}$.
% is a poset in groupoids which is a subcategory of  $\roman{Bl}_{\roman{rs}}$.

\smallskip

(b) For any object $(X_2, I, D_2)$ (cf. Def. 1.4.1 above) of $\Cal{M}$,   the divisor $D_2$
is a simple normal crossings divisor. Sets of labels $S$, together with their functorial behaviour, may be included as separate elements  of the structure of $\Cal{M}$.

\medskip

{\bf 2.2. Two classes of morphisms in $\Cal{M}$.} Fix a category $\Cal{M}$ as above.
 Let $\Cal{X} := (X_2,I_X, D_1)$ and $\Cal{Y} := (Y_2, I_Y, D_2)$ be two objects of $\Cal{M}$.
 
 \smallskip
 
 Assume that we have a locally closed embedding $Y_2\hookrightarrow X_2$
 which extends in a natural way to a morphism between some blowings up of
 $X$, resp $Y$.
 The resulting commutative diagrams will be declared some new morphisms in $\Cal{M}$,
 ``morphisms of closed embeddings'': cf. Sec. 1.5 of [CMMar20].
 
 \smallskip
 
 Similarly, ``morphisms of complements to locally closed embeddings'' are extensions 
 of this definition to $X_2\setminus \overline{Y}_2 \hookrightarrow X_2$ where
$\overline{Y}_2$ denotes the closure of $Y_2$ in $X_2$.
 
 \smallskip
 
 They are presented below as left and right sides of the commutative diagram:
$$
Y_1 \quad \hookrightarrow X_1 \quad\quad  \hookleftarrow X_1\setminus \overline{Y}_1
$$
$$
\downarrow \varphi_Y \ \quad \downarrow \varphi_X \quad\quad \ \downarrow  \ \varphi_{X\setminus Y}
$$
$$
Y_2 \quad \hookrightarrow X_2 \quad\quad  \hookleftarrow X_2\setminus \overline{Y}_2
$$
Subschemes $Y_i\hookrightarrow X_i$  (resp. $X_i\setminus \overline{Y}_i \hookrightarrow X_i$), $i=1,2$,
will be called {\it locally closed} (resp. {\it open}) strata of $X_i$.

\medskip

{\bf 2.3. Example: Kapranov's presentation of $\overline{M}_{0,n}$, $n\ge 3$.} This presentation
of $\overline{M}_{0,n}$, $n\ge 3$, as a result of successive blowings of projective subspaces
in $\bold{P}^{n-3}$ was given in [Ka93], and then used in [BrMe13] for calculating
of regular automorphisms of these stable modular spaces.

\medskip

{\bf 2.4. Example: Connes--Kreimer Hopf algebras from rooted trees.} The Connes--Kreimer
Hopf algebras were introduced in [CoKr00]. 

\smallskip

Later, using the operadic formalism,  F. Chapoton and M. Livernet have  shown
that they appear as well in a geometric environment as above: see
iSec. 6 of [ChaLiv07].

\bigskip

\centerline{\bf 3. NORI MOTIVES AND BURNSIDE GROUPS}

\medskip

Below we work over  a fixed field  $k$ of characteristic 0.
\medskip

{\bf 3.1. Definition.} ([KoTsch19], Sec. 4, Def. 10). {\it Let $\Cal{B}$ (``base scheme'')
be a separated scheme of finite type over $k$.

\smallskip

Consider a smooth $\Cal{B}$--scheme $f: X \to \Cal{B}$. If $U \hookrightarrow X$
is an open embedding with $\overline{U} =X$, then $f|_{U} :U\to \Cal{B}$
is also a smooth $\Cal{B}$--scheme.

\smallskip

(a) Define the set $\roman{Burn_{+}}(\Cal{B})$ as the set of equivalence
classes of  smooth $\Cal{B}$--schemes modulo equivalence relation
generated by $f \sim f|_{U}$ as above. We may denote the respective
equivalence class by $[f: X\to \Cal{B}]$, or simply $[f]$.

\smallskip

(b) Define the monoid structure $+$ upon  $\roman{Burn_{+}}(\Cal{B})$
as generated by disjoint union of smooth $\Cal{B}$--schemes.

It generates the respective Grothendieck group $\roman{Burn} (\Cal{B})$.

\smallskip

(c) Both Burnside group and Burnside monoid a naturally graded: class of 
$X$ of pure dimension $n$ belongs to  $\roman{Burn_{+,n}}(\Cal{B})$
and $\roman{Burn}_n (\Cal{B})$}

\medskip

These constructions are covariant functors of $\Cal{B}$:
a morphism $g:\Cal{B}^{\prime} \to \Cal{B}$ induces
maps $g_*[f] := [g\circ f]$.

\medskip

{\bf 3.2. Boundary homomorphisms.} Start with a pair $Z\subset X$,
in which $X$ is an equidimensional algebraic variety, and $Z$ its closed
subvariety of strictly lesser dimension. Moreover, we will assume that $X$
is reduced and separated, but nothing more. 

\medskip

{\bf 3.2.1. Theorem.} ([KoTsch19], Sec. 4, Theorem 11). {\it On the set 
of Burnside groups of members of such pairs
$Z\subset X$ one can define graded boundary elements

$$
\partial_Z(X) \in \roman{Burn}_{dim (X) -1} (Z)
$$
satisfying two requirements:

\smallskip

(a) For any proper surjective morphism $g: X^{\prime}\to X$ inducing birational equivalence
between $X$ and $X^{\prime}$ and such that $Z^{\prime} = g^{-1}(Z)$,
we have
$$
\partial_Z(X) = (g|_Z)_{*})(\partial_{Z^{\prime}}(X^{\prime})).
$$

(b) If $X$ is smooth, and $Z$ is an $S$--labelled simple normal crossings divisor
in the sense explained in Sec. 1.4 above, then we have an explicit presentation
$$
\partial_Z(X) = - \sum_{\emptyset \ne T\subseteq S} (-1)^{card\,T} [f_T]
$$
where $f_T: D_T \times \bold{A}^{card\, T - 1} \to Z$ is the composition
of projection to $D_T$ and its inclusion into $Z$.

\smallskip

Moreover, these two requirements uniquely determine boundary elements.
}
\smallskip

We can now state precisely and prove Theorem 0.2.1.

\smallskip
Consider a category $\Cal{M}$ as in Sec. 2.1 above,
with its  objects  graded by dimension.
For each object $\Cal{B}$ of $\Cal{M}$, construct its Grothendieck--Burnside group
$\roman{Burn} (\Cal{B})$. Denote by $\roman{Gr Ab}$ the category of graded abelian groups.

\medskip

{\bf 3.3. Theorem.} {\it The natural degree zero map $\roman{Ob} (\Cal{M}) \to \roman{Ob}\,(\roman{Gr Ab})$:  $\Cal{B} \mapsto
\roman{Burn} (\Cal{B})$  extends to the homological Nori geometric diagram, in which
boundary edges correspond to three step towers of closed embeddings $Z\subset Y\subset X$.
The morphism
$$
(X,Y, i+1) \to (Y, Z, i)
$$
sends boundary element $\partial_Y(X)$ to $\partial_Z(Y)$.
}
\medskip

See [MaMar18], subsections 4.7--4.8.

\bigskip

\centerline{\bf 4. BOST--CONNES SYSTEMS AND BURNSIDE GROUPS}

\medskip

{\bf 4.1 Assemblers.} I. Zakharevich introduced and developed in [Za14], [Za15]
a categorical formalism useful for our study.
\smallskip

Briefly, an assembler $\Cal{C}$ is a small category endowed with Grothendieck topology
and an initial object $\emptyset$
satisfying the following restrictions:

\smallskip

(a) All morphisms in it are monomorphisms.

\smallskip

(b) For each object $X$ in $\Cal{C}$ and any two finite disjoint covering families
of $X$ have a common refinement which is also a disjoint covering family.

\smallskip

Here two morphisms $f : Y \to X$ and $g : Z \to X$ are called disjoint
if $Y \times_X  Z = \emptyset $.

\smallskip

An {\it epimorphic assembler $\Cal{C}$ with a sink object $S$} is an assembler,
such that
\smallskip

(c) For each object $X$ of $\Cal{C}$, $Hom_{\Cal{C}}(X,S)$ is non--empty.

\smallskip

(d) Each morphism $f : X\to Y$ in $\Cal{C}$ with non--initial $X$ is an
epimorphism, and this epimorphism is
a covering family.

\smallskip

(e) If $Y,Z \neq \emptyset$, then no two morphisms $Y\to X$ and $Z\to X$
are disjoint.

\medskip

In [MaMar18] it was shown that there is an assembler category,
in the sense of [Za14], underlying the Kontsevich--Tschinkel 
Burnside group, and its equivariant version. Moreover, there is a
Bost--Connes system of endomorphisms acting on the $\widehat{\bold{Z}}$--equivariant
version of the Kontsevich--Tschinkel Burnside ring. This is a lift of the
integral Bost--Connes algebra of [CoConMar08], through the map to the
graded ring associated to the filtration of the Grothendieck ring by
dimension, which in turn maps to the integral Bost--Connes algebra
through the equivariant Euler characteristic. Moreover, this Bost--Connes
structure  can be lifted higher, to endofunctors of the assembler category and
endomorphisms of the associated homotopy--theoretic spectrum.

\smallskip

On the other hand, in [LiMaMar19] it was also shown that the lift of the
Bost--Connes algebra to the $\widehat{\bold{Z}}$--equivariant Grothendieck
ring of varieties can be lifted higher, to a $\widehat{\bold{Z}}$--equivariant category
of Nori motives, where it maps through the fibre functor to the
categorification of the Bost--Connes system constructed in [MarTa14].

\smallskip

We review here briefly the setting of [LiMaMar19] and [MaMar18],
and then we show that a similar categorification of the Bost--Connes
system via Nori motives can be constructed in the case of the 
Kontsevich--Tschinkel Burnside group. As was done in [LieMaMar19] for the
Grothendieck ring, we work here with the relative version of the
Kontsevich--Tschinkel Burnside ring considered in the previous
sections, and with its $\widehat{\bold{Z}}$--equivariant version.

\smallskip

For  much more detailed expositions, see [Za14], [Za15], and Sec. 4 of 
[MaMar18].

\medskip

{\bf 4.2. Relative equivariant Kontsevich--Tschinkel Burnside group.}
In the previous sections we considered a base scheme $\Cal{B}$ and the Kontsevich--Tschinkel
Burnside group $\roman{Burn}(\Cal{B})$.

\smallskip

Here we introduce its equivariant version $\roman{Burn}^{\widehat{\bold{Z}}}(\Cal{B})$
where now $\Cal{B}$ is endowed with a residually finite action of $\widehat{\bold{Z}}$.
The last condition means that this action factors through some finite $\bold{Z}/N\bold{Z}$--quotient of 
$\widehat{\bold{Z}}$. 
\smallskip

 The Burnside group $\roman{Burn}^{\widehat{\bold{Z}}}(\Cal{B})$
then is generated by equivalence classes of smooth $\Cal{B}$--schemes $f: X \to \Cal{B}$
where $X$ is also endowed with
 a residually finite action of $\widehat{\bold{Z}}$, and $f$ is $\widehat{\bold{Z}}$--equivariant.
 Equivalence classes correspond to the equivalence
relation $f \sim f|_U$ where $U\hookrightarrow X$ is a $\widehat{\bold{Z}}$--equivariant dense open embedding. 

\smallskip

As in [LiMaMar19], we will adopt the notation $(\Cal{B},\alpha_{\Cal{B}})$ for the base scheme
endowed with action $\alpha_{\Cal{B}}: \widehat{\bold{Z}} \times \Cal{B} \to \Cal{B}$, and we similarly write
$(X,\alpha_X)$ and $f: (X,\alpha_X)\to (\Cal{B},\alpha_{\Cal{B}})$  to keep track explicitly of the
$\widehat{\bold{Z}}$--actions. 

\smallskip

In [MaMar18] we constructed an assembler category and the associated
homotopy theoretic spectrum underlying the Kontsevich--Tschinkel Burnside group
$\roman{Burn}^{\widehat{\bold{Z}}}(\bold{K})$, where $\bold{K}$ a field. The relative case we consider
here is similar. We give a brief description of the resulting assembler.

\medskip

{\bf 4.3. Proposition.} {\it Consider the family
$\Cal{C}_{(\Cal{B}, \alpha_{\Cal{B}})}^{\widehat{\bold{Z}}}$
of epimorphic assemblers with sink 
whose objects
are $\widehat{\bold{Z}}$--equivariant dense open embeddings $\iota: (U,\alpha_U)\hookrightarrow (X,\alpha_X)$.

Consider the coproduct over equivalence classes of  $\Cal{B}$--schemes 
$f: (X,\alpha_X) \to (\Cal{B},\alpha_{\Cal{B}})$ of members of this family,

The set $\pi_0$ of the associated spectrum is the relative equivariant Kontsevich--Tschinkel set
$\roman{Burn}^{\widehat{\bold{Z}}}(\Cal{B})$.

More generally, the spectra of these assemblers satisfy 
$$
K (\Cal{C}_{(\Cal{B},\alpha_{\Cal{B}} )}^{\widehat{\bold{Z}} }(X,\alpha_X))
\simeq \Sigma_{+}^\infty BG_{(\Cal{B},\alpha_{\Cal{B}} )}^{\widehat{\bold{Z}} }(X,\alpha_X) ,
$$
where 
$G^{\widehat{\bold{Z}}}_{(\Cal{B},\alpha_{\Cal{B}} )} 
(X,\alpha_X)$ is the group of $\widehat{\bold{Z}}$--equivariant birational automorphisms.
}

\smallskip

{\it Proof.} According to Theorem 5.3 of [MarTa14] and Lemma 4.5 of [MaMar18], 
$\Cal{C}^{\widehat{\bold{Z}}}_{(\Cal{B},\alpha_{\Cal{B}} )}(X,\alpha_X)$ is indeed
an epimorphic assembler with sink.

\smallskip

Therefore, in view of Theorem 4.8 of [Za14], we have
$$
\Cal{C}^{\widehat{\bold{Z}}}_{(\Cal{B},\alpha_{\Cal{B}} )}(X,\alpha_X) \simeq
\Sigma_{+}^\infty BG^{\widehat{\bold{Z}}}_{(\Cal{B},\alpha_{\Cal{B}} )} 
(X,\alpha_X),
$$
where  $G^{\widehat{\bold{Z}}}_{(\Cal{B},\alpha_{\Cal{B}} )} (X,\alpha_X)$
is the group of $\widehat{\bold{Z}}$--equivariant birational automorphisms of
the $\Cal{B}$--scheme $f: (X,\alpha_X) \to (\Cal{B},\alpha_{\Cal{B}})$.

\smallskip

The remaining statements also follow from Lemma 4.5 of [MaMar18].

\bigskip

{\bf 4.4. Lifting the Bost--Connes system.} We will use notation and conventions of
[LieMaMar19], which are briefly repeated below.

\smallskip

Put $Z_n:= \roman{Spec}\,K^n$ where $K$ is our ground field.

\smallskip

Given $(\Cal{B},\alpha_{\Cal{B}})$ as above, denote by $\Phi (\alpha_{\Cal{B}})$
the $\widehat{\bold{Z}}$--action upon $\Cal{B}\times Z_n$through composition with
geometric {\it Verschiebung} map: see [LieMaMar19], (2.12) and (2.13).

\smallskip

This action is residually finite as well: if the action $\alpha_{\Cal{B}}$ of 
$\widehat{\bold{Z}}$ factors through $\bold{Z}/N\bold{Z}$, 
then the action $\Phi (\alpha_{\Cal{B}})$ 
 factors through $\bold{Z}/Nn\bold{Z}$. Denote by $\sigma_n$ 
 the endomorphism of multiplication by $n$.
 
 \smallskip
 
 Now, by lifting $\sigma_n$ (and avoiding extra notation),
  define the maps $\sigma_n$ and $\tilde{\rho}_n$ as follows:
$$
\sigma_n :\  (f: (X,\alpha_X) \to (\Cal{B},\alpha_{\Cal{B}})) \mapsto
(f: (X,\alpha_X\circ \sigma_n) \to (\Cal{B},\alpha_{\Cal{B}} \circ \sigma_n)) ,
$$
$$
\tilde{\rho}_n :\  (f: (X,\alpha_X) \to (\Cal{B},\alpha_{\Cal{B}})) \mapsto
  (f \times id : (X\times Z_n, \Phi_n (\alpha_X) \to (\Cal{B} \times Z_n, \Phi_n (\alpha_{\Cal{B}}))  .
$$

{\bf 4.4.1. Proposition.} {\it The maps above can be lifted to functors
of respective assembler categories with sinks that we will denote
by the same letters $\sigma_n$ and $ \tilde{{\rho}}_n$.

These functors induce ring homomorphisms
$$
\sigma_n : \roman{Burn}^{\widehat{\bold{Z}}}(\Cal{B}, \alpha_{\Cal{B}}) \to
\roman{Burn}^{\widehat{\bold{Z}}}(\Cal{B}, \alpha_{\Cal{B}}\circ  \sigma_n)
$$
and group homomorphisms
$$
\tilde{\rho}_n : \roman{Burn}^{\widehat{\bold{Z}}}(\Cal{B}, \alpha_{\Cal{B}}) \to
\roman{Burn}^{\widehat{\bold{Z}}}(\Cal{B} \times Z_n, \Phi_n (\alpha_{\Cal{B}})
)
$$
by the induced maps on the $\pi_0$ of the associated spectra.
}
\smallskip

For proofs, see [MaMar18], Proposition 4.7, and {LieMaMar19], Theorem 3.15.
\smallskip

The functors $\sigma_n$ and $\tilde\rho_n$ of Proposition 4.4.1 and the induced
morphisms of the Kontsevich--Tschinkel Burnside set lift the maps $\sigma_n$
and $\tilde\rho_n$ of the integral Bost--Connes algebra of  [ConCoMar08] in the following way.

\medskip

{\bf Proposition 4.5.} {\it
The action of the $\sigma_n$ and $\tilde{\rho}_n$ on the ${\roman{Burn}}^{\hat{\Z}}(\Cal{B},\alpha_\Cal{B})$
is compatible with the $\sigma_n$ and $\tilde{\rho}_n$ constructed in $\roman{[LieMaMar19]}$
on the Grothendieck rings $K_0^{\hat\Z}(\Cal{B},\alpha_\Cal{B})$, through the 
morphism ${\roman{Burn}^{\hat\Z}}(\Cal{B},\alpha_{\Cal{B}})\to {\roman{gr}} K_0^{\hat\Z}(\Cal{B},\alpha_{\Cal{B}})$ to
the associated graded object  with respect to filtration  by dimension.}

\smallskip

{\it Proof.} The maps $\sigma_n$ and $\tilde{\rho}_n$ on $K_0^{\hat\Z}(\Cal{B},\alpha_\Cal{B})$ are constructed in
 [LieMaMar19] using the same formulas as above. These preserve the filtration
 by dimension, since the schemes $Z_n$ are zero dimensional. Thus, the map 
 ${\roman{Burn}^{\hat\Z}}(\Cal{B},\alpha_\Cal{B})\to {\roman{gr}} K_0^{\hat\Z}(\Cal{B},\alpha_\Cal{B})$ intertwines
 the action of the $\sigma_n$ and $\tilde{\rho}_n$ on the Burnside ring and on the Grothendieck ring.

\medskip
It is also shown in [LieMaMar19], Theorem 2.11, that the actions of  $\sigma_n$ and $\tilde{\rho}_n$ on 
$K_0^{\hat\Z}(\Cal{B},\alpha_\Cal{B})$ are compatible with the $\sigma_n$ and $\tilde\rho_n$
acting on $\Z[\Q/\Z]$, through an equivariant Euler characteristic map .

\bigskip

{\bf 4.6. Equivariant Nori motives.} Given a base scheme $\Cal{B}$ with a residually finite $\hat\Z$--action $\alpha_\Cal{B}$, we
can consider, as in [LieMaMar19], Sec. 7, a category $\Cal{M}_{(\Cal{B},\alpha_\Cal{B})}$ of $\hat\Z$--equivariant Nori motivic
sheaves, in the sense of [Ar08]. This is defined as the abelian category constructed 
from a Nori diagram $\Cal{D}(\Cal{B},\alpha_\Cal{B})$ with the following vertices and edges:

\medskip

$\bullet$  Vertices are given by elements of the form 
$$ 
(f: (X,\alpha_X) \to (\Cal{B},\alpha_\Cal{B}), (Y,\alpha_X|_Y), i, w) 
$$ 
where $f: (X,\alpha_X) \to (\Cal{B},\alpha_\Cal{B})$ is
a smooth $\hat\Z$--equivariant $\Cal{B}$-scheme, with a $\hat\Z$--equivariant embedding
$(Y,\alpha_X|_Y) \hookrightarrow (X,\alpha_X)$,  and $i\in \bold{N}$, $w\in \Z$. 

\medskip

$\bullet$ Edges are of three types:

\smallskip

(1) edges $h^*: (X^{\prime}  \to \Cal{B}, Y^{\prime}, i,w) \to (X\to \Cal{B}, Y, i,w)$ associated to morphisms
$h: (X\to \Cal{B}, Y) \to (X^{\prime} \to \Cal{B}, Y^{\prime})$;

(2) connecting morphisms $\partial: (Y\to \Cal{B},Z,i,w) \to (X\to \Cal{B},Y, i+1,w)$ associated to
a chain of embeddings $Z\hookrightarrow Y \hookrightarrow X$;

(3) twisted projections $(X,Y,i,w) \to (X\times \P^1, Y\times \P^1 \cup X \times \{ 0 \}, i+2, w+1)$.

\medskip

As was shown in [LieMaMar19],  the maps $\sigma_n$ and $\tilde\rho_n$ from above 
induce functors of the categories of Nori motivic sheaves considered above, with
$$ 
\sigma_n: \Cal{M}_{(\Cal{B},\alpha_\Cal{B})} \to \Cal{M}_{(\Cal{B},\alpha_\Cal{B}\circ \sigma_n)} 
$$
$$
 \tilde{\rho}_n: \Cal{M}_{(\Cal{B},\alpha_\Cal{B})} \to  \Cal{M}_{(\Cal{B}\times Z_n,\Phi_n(\alpha_\Cal{B}))}. 
 $$
The functors $\sigma_n$ are compatible with the monoidal structure but $\tilde\rho_n$, are not,
as discussed in [LieMaMar19].

\medskip

By construction, the Bost--Connes structure on the categories $\Cal{M}_{(\Cal{B},\alpha_\Cal{B})}$
of Nori motives is compatible with those discussed above on 
$K_0^{\hat\Z}(\Cal{B},\alpha_\Cal{B})$ and on ${\roman{Burn}^{\hat\Z}}(\Cal{B},\alpha_\Cal{B})$. 
In particular, the same construction can also be applied to the $\hat\Z$-equivariant
version of the Nori motives considered in relation to the categories of blowups.  

\medskip

{\bf 4.7. The Bost--Connes structure of the Kontsevich-Pestun-Tschinkel modular symbols}

A summarized above, in [LieMaMar19] and [MaMar18] we considered various lifts of
the Bost--Connes algebra to Grothendieck rings, assembler, spectra, and Nori motives,
based on varieties with good actions of $\hat\Z$. This setting includes the case of
equivariant Kontsevich--Tschinkel Burnside ring and associated assembler and Nori motives
described here in the previous subsections. There is another setting in birational geometry
where the Bost--Connes structure naturally appears, which we discuss in this subsection,
namely the Kontsevich-Pestun-Tschinkel modular symbols of [KPT19].

\smallskip

As in [KPT19] we consider, for $n\in \N$ and $G$ a finite abelian group, 
the $\Z$-modules $\cM_n(G)$, generated by symbols $\langle a_1, \ldots, a_n \rangle$
with $a_i \in A =G^\vee=Hom(G,\C^*)$, such that $a_1,\ldots, a_n$ generate $A$,
with relations:

(1) $\langle a_{\sigma(1)}, \ldots, a_{\sigma(n)} \rangle =\langle a_1, \ldots, a_n \rangle$,
for all permutations $\sigma\in S_n$;

(2) for all $2\leq k \leq n$ and all $a_1,\ldots, a_k$ and $b_1, \ldots, b_{n-k}$ in $A$ satisfying
$$ \sum_i \Z a_i + \sum_j \Z b_j = A $$
one has the relation
$$ \langle a_1, \ldots, a_k, b_1, \ldots, b_{n-k} \rangle = \sum_{1\leq i \leq k} \langle
a_1-a_i, \ldots, a_i, \ldots, a_k -a_i, b_1, \ldots, b_{n-k} \rangle, $$
with $a_i$ in the $i$th place.

As shown in [KPT19] these relations reflect certain scissor-congruence relations on convex cones in lattices. 
We also consider as in [KPT19] the quotient $\cM_n^-(G)$ of $\cM_n(G)$ by the further relation
$\langle -a_1, \ldots, a_n \rangle=-\langle a_1, \ldots, a_n \rangle$.

\smallskip

There is a $\Z$-bilinear multiplication, for $n=n'+n''$, 
$$ \nabla: \cM_{n'}(G') \otimes \cM_{n''}(G'') \to \cM_n(G) $$
associated to an exact sequence of finite abelian groups
$$ 0 \to G' \to G \to G'' \to 0, $$
with both $G'$ and $G''$ nontrivial, given by
$$ \nabla: \langle a_1, \ldots, a_{n'} \rangle \otimes \langle b_1, \ldots, b_{n''} \rangle \mapsto \sum
\langle \tilde a_1, \ldots, \tilde a_{n'},  \tilde b_1, \ldots, \tilde b_{n''} \rangle, $$
where the sum is over all the lifts $\tilde a_i$  in $A$, in the dual exact sequence
$$ 0 \to A'' \to A \to A' \to 0, $$
while the $\tilde b_j$ are the images in $A$ of the $b_j$ under the embedding $A''\to A$. 
The multiplication map descends to the $\cM_n^-(G)$ in the same form.

\smallskip

There is similarly a $\Z$-bilinear comultiplication, for $n=n'+n''$, determined by
$$ \Delta: \cM_n(G) \to \cM_{n'}(G') \otimes \cM^-_{n''}(G'') $$
for a sequence as above with $G''$ nontrivial, given by
$$ \Delta: \langle a_1, \ldots, a_n \rangle \mapsto \sum \langle a_{I'} \mod A'' \rangle \otimes \langle a_{I''} \rangle^-, $$
for $\{ 1, \ldots, n\} =I' \sqcup I''$ with $\# I'=n'$ and $\# I''=n''$, and $a_{I'}=a_{i_1},\ldots, a_{i_{n'}}$
and $a_{I''}=a_{j_1}\ldots a_{j_{n''}}$ for $I'=\{ i_1,\ldots, i_{n'}\}$ and $I''=\{ j_1,\ldots, j_{n''}\}$,
such that all the $a_{j_k}$ are in $A''\hookrightarrow A$.
Here taking the quotient $\cM^-_{n''}(G'')$ instead of $\cM_{n''}(G'')$ is necessary because of the
second type of relations in $\cM_n(G)$, see Proposition~9 of [KPT19].
The comultiplication also descends to $\cM_n^-(G)$.

\smallskip

In particular, we will focus here on the case of $\cM_{n,N}:=\cM_n(\Z/N\Z)$. We write $G_N=\Z/N\Z$
and $A_N$ for the characters and we consider the projective system of the $A_N$ ordered by
divisibility, with the maps $\sigma_k: A_N \to A_M$ for $M|N$ with $N=Mk$, given by
$\sigma_k : \zeta \mapsto \zeta^k$ when we identify $A_N$ with the group of 
$N$th roots of unity (multiplication by $k$ if written additively).  Dually we have the injective 
system of the $G_N$ with the corresponding inclusions $j_k: G_M \to G_N$, so that
$$ \varinjlim G_N =\Q/\Z, \ \ \ \ \   \varprojlim A_N =\hat\Z. $$

\smallskip

We define $\cM_n(\Q/\Z)$ as the $\Z$-module spanned by the $\cM_n(G_N)$ for
all $N$. We write elements of $\cM_n(\Q/\Z)$ in the form $\sum_i c_i \langle a_i \rangle$
where we use the shorthand notation $\langle a_i \rangle:=\langle a_{i,1}, \ldots, a_{i,n} \rangle$ 
with $a_{i,k}\in A_{N_i}$. We also write $\cM_\Q(\Q/\Z):=\cM(\Q/\Z)\otimes_\Z \Q$.
We also write $\cM_n(\Q/\Z)^-$ for the span of the $\cM_n(G_N)^-$.

\smallskip

Let $\sigma_k: \cM_n(\Q/\Z) \to \cM_n(\Q/\Z)$ be the map of $\Z$-modules
determined by
$$ \sigma_k: \langle a \rangle \mapsto \langle \sigma_k(a) \rangle, $$
for $\langle a \rangle =\langle a_1, \ldots, a_n \rangle$ in $\cM_n(\Z/N\Z)$ and
$\langle \sigma_k(a) \rangle=\langle \sigma_k(a_1), \ldots, \sigma_k(a_n) \rangle$.
We set $\sigma_k \langle a \rangle=0$ whenever $\langle \sigma_k(a) \rangle$ would
not be an acceptable symbol: for instance this happens for $\langle a \rangle\in \cM_n(\Z/k\Z)$
when all $\sigma_k(a_i)=0 \in \Z/k\Z$. We also consider the maps of $\Z$-modules 
$$ \rho_k:  \cM_n(G_M) \to \cM_n(G_N) $$
that maps a symbol to the sum over preimages
$$  \rho_k:
\langle a_1, \ldots, a_n \rangle \mapsto \sum_{\sigma_k(b_i)=a_i} \langle b_1, \ldots, b_n \rangle . $$
With a slight abuse of notation we will write for convenience
$$ \langle \rho_k(a) \rangle:= \sum_{\sigma_k(b_i)=a_i} \langle b_1, \ldots, b_n \rangle, $$
so that we can write $\tilde\rho_k \langle a \rangle=\langle \tilde\rho_k(a) \rangle$.
On $\cM_\Q(\Q/\Z)$ we also consider the $\Q$-linear maps
$$ \hat\rho_k =\frac{1}{k^n} \tilde\rho_k:  \cM_{\Q,n}(\Q/\Z) \to \cM_{\Q,n}(\Q/\Z), $$ 
$$ \hat\rho_k : \langle a_1, \ldots, a_n \rangle \mapsto \frac{1}{k^n}
 \sum_{\sigma_k(b_i)=a_i} \langle b_1, \ldots, b_n \rangle. $$
 
 \smallskip
 
 Notice that the notation we are using here differs from [ConCoMar08], [LieMaMar19] , [MaMar18]
 where we used $\tilde\rho_k$ for what is here called $\rho_k$ and $\rho_k$ for what
 is here called $\hat\rho_k$. The reason for this change of notation is that, as maps
 on the group ring $\Z[\Q/\Z]$, the $\tilde\rho_k$ of \cite{CCM} are not ring homomorphisms, 
 while the $\rho_k$ of \cite{CCM}  are algebra homomorphisms of $\Q[\Q/\Z]$ and are the
 original Bost--Connes maps. Here however, since we will be considering the multiplication
 structure on the $\cM_n(G)$ as in [KPT19] rather than the group ring $\Z[\Q/\Z]$, it
 will turn out that the $\tilde\rho_k$ of [ConCoMar08]  are the fundamental maps, which we
 call simply $\rho_k$.

\smallskip

{\bf Lemma 4.8.} {\it
The maps $\sigma_k$ and $\rho_k$ on $\cM_n(\Q/\Z)$ defined as above satisfy 
the relations $\sigma_k \sigma_\ell =\sigma_{k \ell}$, $\rho_k \rho_\ell =\rho_{k\ell}$.
Moreover, if $(k,\ell)=\gcd \{ k, \ell \}=1$, then
$$ \sigma_k \circ \rho_\ell = \rho_\ell \circ \sigma_k,  $$
while 
$$ \sigma_k \circ \rho_k =k\cdot id, \ \ \ \ \rho_k \circ \sigma_k = e_k, $$
with the map $e_k: \cM_n(\Q/\Z)\to \cM_n(\Q/\Z)$ given by
$$ e_k \langle a_1,\ldots, a_n \rangle = \sum_{k s_i =0} \langle a_1+s_1,\ldots, a_n+s_n \rangle. $$ }

{\it Proof.} The compositions $\sigma_k \sigma_\ell =\sigma_{k \ell}$ and $\rho_k \rho_\ell =\rho_{k\ell}$
follow directly from the definition.
For $(k,\ell)=1$ we have $$\sigma_k \sum_{\ell b_i =a_i} \langle b_1, \ldots, b_n \rangle =
\sum_{\ell b_i =a_i} \langle k b_1, \ldots, k b_n \rangle $$ $$ = \sum_{\ell k b_i=k a_i} \langle k b_1, \ldots, k b_n \rangle=
\rho_\ell \langle k a_1, \ldots, k a_n \rangle.$$ On the other hand $\sigma_k \sum_{k b_i =a_i} \langle b_1, \ldots, 
b_n \rangle = \sum_{k b_i =a_i} \langle a_1, \ldots, a_n \rangle = k \cdot \langle a_1, \ldots, a_n \rangle$ and
$\rho_k \langle k a_1, \ldots, k a_n \rangle=\sum_{k b_i=k a_i} \langle b_1, \ldots, b_n \rangle=
\sum_{k s_i=0} \langle a_1+s_1, \ldots, a_n +s_n\rangle$.  The more general case with 
$(k,\ell)\neq 1$ is also obtained by combining the relations above.

\smallskip

Note that, unlike in the case of the maps $\sigma_n$ and $\rho_n$ acting on $\Q[\Q/\Z]$ in
the Bost--Connes system, the $e_k$ here are not idempotents with respect to $\nabla$.

\smallskip

The following lemma is a direct consequence of the explicit expression for the
multiplication and comultiplication and of the maps $\sigma_k$ and $\rho_k$
defined above.

\smallskip

{\bf Lemma 4.9.} {\it
The multiplication $\nabla$ induces a multiplication 
$$ \nabla: \cM_{n'}(\Q/\Z)\otimes \cM_{n''}(\Q/\Z) \to \cM_{n'+n''}(\Q/\Z), $$
which is given by
$$ \nabla (\langle a \rangle, \langle b \rangle) = \langle \rho_\ell(a), b \rangle, $$
for $\langle b \rangle \in \cM_n(\Z/\ell\Z)$. The comultiplication $\Delta$ induces
a comultiplication 
$$ \Delta: \cM_n(\Q/\Z) \to \oplus_{n'+n''=n} \cM_{n'}(\Q/\Z) \otimes \cM_{n''}(\Q/\Z) $$
given by
$$ \Delta \langle a \rangle = \oplus \langle \sigma_\ell(a_{I'}) \rangle \otimes \langle a_{I''} \rangle^-, $$
where for $\langle a \rangle\in \cM_n(G_N)$ the sum is over $\ell | N$ and over all partitions $$I'\sqcup I'' =\{ 1, \ldots, n\},$$
with $\langle a_{I''} \rangle$ in $\cM_n(Ker(\sigma_\ell))$. }

\smallskip

Thus, we can consider the graded ring $\cM(\Q/\Z)=\oplus_n \cM_n(\Q/\Z)$ endowed with this
multiplication, and the $\Q$-algebra $\cM_\Q(\Q/\Z)$. 

\smallskip

{\bf Proposition 4.10.} {\it 
The maps $\rho_k : \cM(\Q/\Z) \to \cM(\Q/\Z)$ are ring homomorphisms, while 
the maps $\sigma_k : \cM(\Q/\Z) \to \cM(\Q/\Z)$ are coalgebra homomorphisms.}

\smallskip

{\it Proof.}  As mentioned above, the multiplication in $\cM(\Q/\Z)$ is given by 
$\nabla (\langle a \rangle, \langle b \rangle) = \langle \rho_\ell(a), b \rangle$, hence we have
$$ \rho_k \circ \nabla (\langle a \rangle, \langle b \rangle) = \langle \rho_k \rho_\ell (a), \rho_k(b) \rangle=
\langle \rho_\ell \rho_k (a), \rho_k(b) \rangle=\nabla \circ \rho_k (\langle a \rangle, \langle b \rangle), $$
while the comultiplication is given by 
$\Delta \langle a \rangle = \oplus \langle \sigma_\ell(a_{I'}) \rangle \otimes \langle a_{I''} \rangle^-$ so that
$$ \sigma_k \circ \Delta \langle a \rangle = \oplus \langle \sigma_k \sigma_\ell(a_{I'}) \rangle \otimes \langle \sigma_k (a_{I''}) \rangle^- . $$
The sum is over all $\ell | N$, for $\langle a \rangle\in \cM_n(\Z/N\Z)$ and over all decompositions $I',I''$
where $\langle a_{I''} \rangle$ is in $\cM_n(\Z/\ell\Z)$. If $(k,N)=1$, then also $(k,\ell)=1$ and we have 
$\sigma_\ell : \Z/N\Z \to \Z/N\Z$ and $\langle a_{I''} \rangle$ is in $\cM_n(\Z/\ell\Z)$ iff the same holds for
$\langle \sigma_k (a_{I''}) \rangle$, hence the sum on the right hand side is the same as
$\Delta \circ \sigma_k \langle a \rangle$. It suffices then to check the case where $k|N$, with $N=kM$.
Then $\langle \sigma_k(a) \rangle \in \cM_n(\Z/M\Z)$ and $\Delta \langle \sigma_k(a) \rangle$ is the
sum over all the $\ell'|M$ and the partitions $I',I''$ with $\langle b_{I''} \rangle \in \cM_n(\Z/\ell'\Z)$,
up to terms that are mapped to zero by $\sigma_k$ when $a_{I''}\in \Z/k\Z$.
Thus, the maps $\rho_k$ are compatible with the multiplication $\nabla$ and the maps $\sigma_k$
are compatible with the comultiplication $\Delta$. 

\smallskip

Note that the maps $\sigma_k$ are not compatible with multiplication since
$$ \sigma_k \nabla (\langle a \rangle, \langle b \rangle) = \langle \sigma_k \rho_\ell (a), \sigma_k(b) \rangle. $$
For $(k,\ell)=1$ one has $\sigma_k \rho_\ell = \rho_\ell \sigma_k$, which gives the right compatibility,
but in the case where $(k,\ell)\neq 1$ this is not the case. For simplicity, consider the case where $\ell=k$,
where $\langle \sigma_k \rho_k (a), \sigma_k(b) \rangle=\langle a, \sigma_k(b) \rangle$ while
$\langle \rho_k \sigma_k (a), \sigma_k(b) \rangle =\sum_{k s_i =0} \langle a+s, \sigma_k(b) \rangle$,
with $a+s=(a_1+s_1,\ldots,a_n+s_n)$. Similarly, the maps $\rho_k$ are not compatible with
comultiplication for exactly the same reason.

\smallskip

The data $(\cM(\Q/\Z),\nabla,\Delta,\rho_k,\sigma_k)$ provide the Bost--Connes structure associated
to the Kontsevich-Pestun-Tschinkel modular symbols. 

\smallskip

In the original Bost--Connes quantum statistical mechanical system, the maps
$\rho_k$ have a geometric interpretation in terms of lattices and commensurability,
see Chapter~3 of  [CoMa08]. One defines $\Q$-lattices as pairs $(L,\phi)$ 
of a rank $n$-lattice $L$ (with $n=1$ in the original Bost--Connes case) and a
group homomorphism $\phi: \Q^n/\Z^n \to \Q L/L$. Two $\Q$-lattices are
commensurable if $\Q L_1=\Q L_2$ and $\phi_1=\phi_2$ modulo $L_1+L_2$.
In the rank one case, the commensurability relation is precisely captured by
the maps $\rho_k$ of the Bost--Connes system. 

\smallskip

In the case of the Kontsevich-Pestun-Tschinkel 
modular symbols, it is shown in [KPT19] that the $\cM_n(G)$ can also be
described in terms of lattices and cones.
For $\cM_n(G)$, one considers isomorphism classes of data $(L, \chi, \Lambda)$ with
$L\simeq \Z^n$ a lattice of rank $n$, with an element $\chi\in L\otimes A$ for which
the induced homomorphism $L^\vee \to A$ is surjective, and a stricly convex cone 
$\Lambda \subset L_\R$ spanned by a basis of $L$. To a datum $(L,\chi,\Lambda)$
one associates a symbol $\psi(L,\chi,\Lambda)=\langle a_1,\ldots, a_n \rangle$ obtained
by a choice of a basis $\{ e_i \}$ of $L$ with respect to which $\chi=\sum_i e_i \otimes a_i$.
The relations of $\cM_n(G)$ are then described geometrically in terms of scissor congruence
relations $\psi(L,\chi,\Lambda)=\sum_i \psi(L,\chi,\Lambda_i)$ determined by possible
decompositions of the cone. 

\smallskip

To see the compatibility between these two settings, note that a lattice $L$ with
the choice of a basis $\{ e_i \}$ determines an invertible $\Q$-lattice of rank $n$,
by taking $\phi : \Q^n/\Z^n \to \Q L/L$ that maps $\phi: (a_1,\ldots, a_n)\mapsto \sum_i a_i e_i$.
To a (not necessarily invertible) $\Q$-lattice with $\phi : \Q^n/\Z^n \to \Q L/L$,
we can assign an element $\chi\in L\otimes \Q/\Z$ by $\chi= \sum_i \phi(e_i) \otimes a_i$
with $e_i$ the tuple in $(\Q/\Z)^n$ with $1$ in the $i$th position and zeroes elsewhere. 
For a fixed $\Z/N\Z$, the restriction of $\phi$ to $(\Z/N\Z)^n$ induces on those
$n$-tuples that determine corresponding symbols $\langle a_1, \ldots, a_n \rangle$
a map $\langle a_1, \ldots, a_n \rangle \mapsto \sum_i e_i \otimes a_i$ that 
gives rise to an associated $\chi$ as above. The maps $\rho_k$ are then implementing
the commensurability relation as in the Bost--Connes case on the rank one $\Q$-lattices
obtained by projecting along the $e_i$ directions.

\bigskip

\newpage
\centerline{\bf References}

\medskip

%[Ar08] D. Arapura.  {\it An abelian category of motivic sheaves.}, Adv. Math. 233 (2013), pp. 135--195. 
%arXiv:0801.0261.

%\smallskip

[AT16] D. Abramovich, M. Temkin. {\it Functorial factorization of birational maps
for qe schemes in characteristic 0.} Journ. of Algebra and Number Theory, 13(2) (2019),
pp. 379--424. arXiv:1606.08424, 45 pp.

\smallskip

[AVi02] D. Abramovich, A. Vistoli. {\it Compactifying the space of stable maps.}
Journ. Amer. Math. Soc., 15(1) (2002), pp. 27--75.

\smallskip

[BeRy19] D.  Bergh, D.  Rydh. {\it Functorial destackification and weak factorization of orbifolds.}
arXiv:1905.00872. 37 pp.

\smallskip

[BoCon95] J.  B.  Bost,  A.  Connes.  {\it Hecke algebras, type III factors and phase transitions with
spontaneous symmetry breaking in number theory}. Selecta Math. (new ser.),
1 (1995), no. 3,  pp. 411--457.

\smallskip

[BrMe13] A. Bruno, M. Mella. {\it The automorphism group of  $\overline{M}_{0,n}$.}
Journ. European Math. Soc. 15 (2013), pp. 959--968.

\smallskip

%[Ca14]  P.  Cartier. {\it A country of which nothing is known but the name:
%Grothendieck and ``motives''.}   In:  Alexandre Grothendieck: A Mathematical Portrait.
%Ed. by Leila Schneps. International Press of Boston (2014),
%pp. 269--298. %

%\smallskip

[ChaLiv07] F.~Chapoton, M.~Livernet. {\it Relating two Hopf algebras built from an operad}.
Int. Math. Res. Not. IMRN, no. 24 (2007),  27 pp. arXiv:0707.3725.

\smallskip

[CMa19] N.  Combe,  Yu.  Manin.  {\it Symmetries of genus zero modular operad.}
To be published in the Memorial Collection in honour of B. Dubrovin.
 arXiv:1907.10317. 11 pp.
 \smallskip
 
[CMaMar20] N. Combe, Yu. Manin,  M. Marcolli. {\it Dessins for modular operad
and Grothendieck--Teichm\"uller  group.} To be published in the Collection in honor of
V. Turaev. arXiv:2006.13663. 28 pp.

\smallskip

[ConCoMar08] A.~Connes, C.~Consani, M.~Marcolli, {\it Fun with $\F_1$}. J. Number Theory 
129 (2009), pp. 1532--1561. arXiv:0806.2401.

\smallskip

[CoKr00]  A.~Connes, D.~Kreimer. {\it Renormalization in quantum field theory and the 
Riemann--Hilbert problem I:
The Hopf algebra structure of graphs and the main theorem}. Comm. Math. Phys., Vol.~210 (2000) N.1, 249--273.

\smallskip

[CoMa08] A.~Connes, M.~Marcolli, {\it Noncommutative geometry, quantum fields and motives},
Colloquium Publications, Vol.55, American Mathematical Society, 2008.

\smallskip

[HuM-St17] A.~Huber, St.~M\"uller--Stach. {\it Periods and Nori motives. With
contributions by Benjamin Friedrich and Jonas von Wangenheim.} Springer, 2017,
xxiii+372 pp.  

\smallskip

[Ka93]  M. Kapranov. {\it Veronese curves and Grothendieck--Knudsen moduli spaces 
 $\overline{M}_{0,n}$.} Journ. Alg. Geom. 2 (1993), 239--262.

\smallskip

[KoTsch19] M. Kontsevich,  Yu.  Tschinkel.  {\it Specialization of birational types.}
Inv. Math. 217 (2019), 415--432.

\smallskip

[KPT19] M.~Kontsevich, V.~Pestun, Yu.~Tschinkel, {\it Equivariant birational geometry and
modular symbols}, arXiv:1902.09894v2.

\smallskip

[KrTsch20] A. Kresch, Yu. Tschinkel. {\it Equivariant birational types and Burnside
volume.} arXiv:2007.12538. 41 pp.

\smallskip

[LieMaMar19] J.F. Lieber, Yu.I. Manin, M. Marcolli. {\it Bost--Connes systems and F1-structures 
in Grothendieck rings, spectra, and Nori motives}, arXiv:1901.00020.

\smallskip

[MaMar18]  Yu. Manin,  M.  Marcolli.  {\it Homotopy types and geometries below  $Spec\, \bold{Z}$.}
In: Dynamics: Topology and Numbers. Conference to the memory of Sergiy Kolyada.
Contemp. Math., AMS 744 (2020), pp. 27--54. 

arXiv:math/CT/1806.10801.

\smallskip

[MaMar19--1]  Yu. Manin,  M.  Marcolli.  {\it Nori diagrams and persistence homology.}
Math. Comput. Sci. 14 (2020), no. 1, pp. 77--102.
arXiv:1901.1031.

\smallskip
[MaMar19--2]  Yu. Manin,  M.  Marcolli.  {\it Quantum statistical mechanics of the
absolute Galois group.} In: Special Issue of Symmetry, Integrability and Geometry: Methods
and Applications, in honor of Motohiko Mulase for his 65th birthday.
SIGMA 16 (2020), 52 pp.    arXiv:mathAG/1907.13545v2. 62 pp. 
\smallskip

[MarTa14] M.~Marcolli, G.~Tabuada. {\it Bost-Connes systems, categorification, quantum statistical mechanics, and Weil numbers}.
 J. Noncommut. Geom. 11 (2017) no. 1, pp. 1--49. 
 arXiv:1411.3223.
\smallskip
[Za14]  I.~Zakharevich.
{\it The $K$-theory of assemblers}. Adv. Math. 304 (2017), pp. 1176--1218.
arXiv:1401.3712.
\smallskip

[Za15] I.~Zakharevich. {\it On $K_1$ of an assembler.} J. Pure Appl. Algebra, 221 (2017), no. 7, pp. 1867--1898.
arXiv:1506197.

\smallskip

\bigskip

{\bf No\'emie C. Combe, Max-Planck-Institut f\"ur Mathematik in den Naturwissenschaften,
Inselstr. 22, 04103 Leipzig, Germany}

\medskip

{\bf Yuri I. Manin, Max-Planck-Institut f\"ur Mathematik, Vivatsgasse 7, 53111 Bonn, Germany}

\medskip

{\bf Matilde Marcolli, Math. Department, Mail Code 253-37, Caltech, 1200 E.California Blvd.,
 Pasadena, CA 91125, USA}

\enddocument